\documentclass[titlepage,twoside,12pt]{article}
\usepackage{amssymb}
\usepackage{amsfonts}
\begin{document}

{\LARGE \bf Extending Mappings between Posets} \\ \\

{\bf Elem\'{e}r E ~Rosinger} \\ \\
{\small \it Department of Mathematics \\ and Applied Mathematics} \\
{\small \it University of Pretoria} \\
{\small \it Pretoria} \\
{\small \it 0002 South Africa} \\
{\small \it eerosinger@hotmail.com} \\ \\

{\bf Abstract} \\

A variety of possible extensions of mappings between posets to their Dedekind order completion
is presented. One of such extensions has recently been used for solving large classes of
nonlinear systems of partial differential equations with possibly associated initial and/or
boundary value problems.  \\ \\ \\

{\bf 1. The General Setup} \\

Let $( X, \leq )$ and $( Y, \leq )$ be two arbitrary posets and \\

(1.1)~~~ $ \varphi : X ~\longrightarrow~ Y $ \\

any mapping between them. We shall be interested to set up {\it commutative diagrams} \\

\begin{math}
\setlength{\unitlength}{0.1cm}
\thicklines
\begin{picture}(50,31)

\put(0,12){$(1.2)$}
\put(23,25){$X$}
\put(50,29){$\varphi$}
\put(30,26){\vector(1,0){46}}
\put(81,25){$Y$}
\put(24,22){\vector(0,-1){16}}
\put(22,0){$X^{\#}$}
\put(30,1){\vector(1,0){46}}
\put(81,0){$Y^{\#}$}
\put(82,22){\vector(0,-1){16}}
\put(50,-4){$\varphi^\diamondsuit$}

\end{picture}
\end{math} \\ \\

where $X^{\#}$ and $Y^{\#}$ are the Dedekind order completions, [3,2,4], of $X$ and $Y$,
respectively, while the mappings \\

(1.3)~~~ $ \varphi^\diamondsuit : X^{\#} ~\longrightarrow~ Y^{\#} $ \\

are {\it extensions} of the given mapping in (1.1), in view of the commutativity of (1.2). \\

As we shall see, there are many natural ways to obtain extensions (1.3). One such way, see
(A.26) - (A.28) and Proposition A.1 in the Appendix, has recently been used successfully in
order to solve large classes of nonlinear systems of PDEs with possibly associated initial
and/or boundary value problems, [4,1,5-7]. \\

Several other earlier obtained results relating to posets and their Dedekind order completions,
result needed in the sequel, are summarized in the Appendix. \\

In view of the main interest pursued being the solution of large classes of nonlinear systems
of PDEs with possibly associated initial and/or boundary value problems, the sets $X$ and $Y$
are supposed to be infinite, since in the particular case when solving PDEs, they correspond
to spaces of functions on Euclidean domains on which the respective PDEs are defined. \\
Furthermore, for the convenience of the Dedekind order completion method, [3], and without
loss of generality, [3,2,4], we shall assume that the posets $( X, \leq )$ and $( Y, \leq )$
do not have minimum or maximum. Otherwise, these two posets can be arbitrary. \\ \\

{\bf 2. Constructing Extensions} \\

It is quite natural to define the extension (1.3) as follows, see (A.7), (A.8) \\

(2.1)~~~ $ {\cal P} ( X ) \ni A ~\longmapsto~ \varphi^\diamondsuit ( A )
                                     ~=~ ( \varphi ( A ) )^{ul} \in Y^{\#} $ \\

which enjoys the following two advantages :

\begin{itemize}

\item it has a larger domain of definition that required in (1.3), and furthermore

\item it does not make use of the partial order on $X$.

\end{itemize}

This however, is precisely the definition of the mapping $\varphi^{\#}$ in (A.26) - (A.28)
which, a mentioned, was given earlier in [4], and used in solving large classes of nonlinear
systems of PDEs with possibly associated initial and/or boundary value problems, [4,1,5-7]. \\
Consequently, we shall look for other possible extensions (1.3) which may similarly be
natural. \\

Let us start by noting that the desired extended mapping $\varphi^\diamondsuit$ in (1.2), (1.3)
must be such that, given, $A \subseteq X$, in order to obtain the corresponding $\varphi^
\diamondsuit ( A ) \subseteq Y$, one should not use more information than it is in the subset
$\varphi ( A ) \subseteq Y$. This is precisely the reason $\varphi^{\#}$ was defined in the
respective manner in (A.26), (A.27), see also (2.1) above. \\

And then, the way left for alternative definitions of $\varphi^\diamondsuit$ is to try to use
in the definition of $\varphi^\diamondsuit ( A ) \subseteq Y$, with $A \subseteq X$, an amount
of information which may possibly be {\it less} than that contained in $\varphi ( A )
\subseteq Y$. \\

A simplest way to do that is to define \\

(2.2)~~~ $ \widetilde \varphi : {\cal P} ( X ) ~\longrightarrow~ Y^{\#} $ \\

by \\

(2.3)~~~ $ \widetilde \varphi ( A ) ~=~
           \bigcap_{\, a \in A}~ ( \varphi (~ [ a > \, \cap \, A ~) )^{ul},~~~
                                                                  A ~\subseteq~ X  $ \\

This definition can obviously be generalized in the following manner. A mapping \\

(2.4)~~~ $ L : {\cal P} ( X ) ~\longrightarrow~ {\cal P} ( X ) $ \\

is called {\it cofinal}, if and only if \\

(2.5)~~~ $ \begin{array}{l}
                   \forall~~ A \subseteq X ~:~ \\ \\
                   ~~~~*)~~ L ( A ) ~\subseteq~ A \\ \\
                   ~~**)~~ L ( A ) ~~\mbox{is cofinal in}~ A
            \end{array} $ \\

Here we recall that a subset $B \subseteq A$ is {\it cofinal} in $A$, if and only if \\

(2.6)~~~ $ \forall~~ a \in A ~:~ \exists~~ b \in B ~:~ a ~\leq~ b $ \\

And then we can define \\

(2.7)~~~ $ \varphi^L :  {\cal P} ( X ) ~\longrightarrow~ {\cal P} ( X ) $ \\

by \\

(2.8)~~~ $ \varphi^L ( A ) ~=~
             \bigcap_{\, a \in L ( A )}~( \varphi (~ [ a > \, \cap \, A ~) )^{ul},~~~
                                                                  A ~\subseteq~ X  $ \\

This further suggests the following alternative possibility. Given $A \subseteq X$, instead of
the subsets $[ a > \, \cap \, A \subseteq A$, with $a \in A$, or $L ( A ) \subseteq A$, we can
consider arbitrary subsets $B \subseteq A$. \\
However, in defining $\varphi^\diamondsuit ( A )$, one should not lose too much from the
information in $\varphi ( A )$. Thus there should be some {\it restriction} on what kind of
subsets $B \subseteq A$ one is considering. \\
In this regard, and as above, a natural candidate is given by subsets $B \subseteq A$ which
are {\it cofinal} in $A$. And then, we arrive at defining \\

(2.9)~~~ $ \overline{\varphi} : {\cal P} ( X ) ~\longrightarrow~ Y^{\#} $ \\

by \\

(2.10)~~~ $ \overline{\varphi} ( A )
                 ~=~ \bigcap_{\, B ~cofinal~in~ A}~ ( \varphi ( B ) )^{ul} $ \\

We note that, unlike $\widetilde \varphi$ and $\overline{\varphi}$ which are two possible
definitions for $\varphi^\diamondsuit ( A )$, there can in general be {\it infinitely} many
mappings $\varphi^L$, for any given pair of posets $( X, \leq )$ and $( Y, \leq )$. \\

We also note that in view of (A.11), one obtains \\

(2.11)~~~ $ \overline{\varphi} ( A ) ~\cup~ \widetilde \varphi ( A ) ~\cup~
                \varphi^L ( A ) ~\subseteq~ \varphi^{\#} ( A ),~~~ A \subseteq X $ \\ \\

{\bf 3. Relations Among the Extended Mappings $\widetilde \varphi,~ \varphi^L,~
\overline{\varphi}$ and $\varphi^{\#}$} \\

{\bf Proposition 3.1.} \\

(3.1)~~~ $ \widetilde \varphi ~=~ \varphi^L $ \\

for every {\it cofinal} mapping $L$ in (2.4). \\

{\bf Proof} \\

In view of (2.5), we have $L ( A ) \subseteq A$, thus (2.3), (2.8) yield the inclusion
'$\subseteq$' in (3.1). \\
For the converse inclusion '$\supseteq$' in (3.1), we recall that $L ( A )$ is cofinal in $A$,
see (2.5). Hence for every $a \in A$, there exists $a\,' \in L ( A )$, such that $a \leq
a\,'$. Consequently, we have \\

$~~~~~~ [ a\,' > \cap\, A ~\subseteq~ [ a > \cap\, A $ \\

thus \\

$~~~~~~ \varphi (  [ a\,' > \cap\, A ) ~\subseteq~ \varphi ( [ a > \cap\, A ) $ \\

and then (A.11) implies \\

$~~~~~~ ( \varphi (  [ a\,' > \cap\, A ) )^{ul} ~\subseteq~
                                       ( \varphi ( [ a > \cap\, A ) )^{ul} $ \\

and the proof of (3.1) is completed. \\

{\bf Proposition 3.2.} \\

Let $A \subseteq X$ be {\it directed}, then \\

(3.2)~~~ $ \overline{\varphi} ( A ) ~\subseteq~ \widetilde \varphi ( A ) $ \\

Here we recall that $A \subseteq X$ is {\it directed}, if and only if \\

(3.3)~~~ $ \forall~~ a,~ a\,' \in A ~:~
                      \exists~~ a\,'' \in A ~:~ a ~\leq~ a\,'',~~ a\,' ~\leq~ a\,'' $ \\

{\bf Proof} \\

In view of (2.3), let $a \in A$, then $B = [ a >\ \cap\, A$ is cofinal in $A$, since $A$ is
directed. Hence (2.10) gives the inclusion in (3.2). \\

{\bf Proposition 3.3.} \\

If the mapping $\varphi$ in (1.1) is {\it increasing}, then \\

(3.4)~~~ $ \overline{\varphi} ~=~ \varphi^{\#} $ \\

{\bf Proof} \\

We shall show that \\

(3.5)~~~ $ ( \varphi ( A ) )^{ul} ~=~ ( \varphi ( B ) )^{ul} $ \\

for every $B \subseteq A$, with $B$ cofinal in $A$. Indeed, for every $a \in A$, there exists
$b \in B$, such that $a \leq b$. But $\varphi$ is increasing, hence $\varphi ( a ) \leq
\varphi ( b )$, which means \\

$~~~~~~ [ \varphi ( b ) > ~\subseteq~ [ \varphi ( a ) > $ \\

thus in view of (A.2), we obtain \\

$~~~~~~ ( \varphi ( B ) )^u ~\subseteq~ ( \varphi ( A ) )^u $ \\

But $B \subseteq A$ and (A.11) always imply \\

$~~~~~~ ( \varphi ( A ) )^u ~\subseteq~ ( \varphi ( B ) )^u $ \\

Hence in our case (3.5) does indeed hold. And then (3.4) follows from (2.10) and (A.27). \\

{\bf Corollary 3.1} \\

If the mapping $\varphi$ in (1.1) is {\it increasing}, then \\

(3.6)~~~ $ \overline{\varphi} ( A ) ~=~ \widetilde\varphi ( A )
                     ~=~ \varphi^L ( A ) ~=~ \varphi^{\#} ( A ) $ \\

for every {\it directed} $A \subseteq X$. \\ \\

{\bf 4. Extension Diagrams} \\

Let us return now to the initial main problem, namely, to construct extensions (1.2) for
arbitrary mappings (1.1) by using Dedekind order completions. \\

{\bf Theorem 4.1} \\

Let $\varphi$ in (1.1) be an arbitrary mapping, then the following two diagrams are
commutative \\

\begin{math}
\setlength{\unitlength}{0.1cm}
\thicklines
\begin{picture}(50,31)

\put(0,12){$(4.1)$}
\put(23,25){$X \ni x$}
\put(62,29){$\varphi$}
\put(38,26){\vector(1,0){52}}
\put(38,25){\line(0,1){2}}
\put(95,25){$\varphi ( x ) \in Y$}
\put(33,22){\vector(0,-1){16}}
\put(32,22){\line(1,0){2}}
\put(14,0){${\cal P} ( X ) \ni \{\, x \,\}$}
\put(40,1){\vector(1,0){49}}
\put(40,0){\line(0,1){2.5}}
\put(93,0){$< \varphi ( x ) \,] \in Y^{\#}$}
\put(100,22){\vector(0,-1){16}}
\put(99,22){\line(1,0){2}}
\put(52,-4){$\overline{\varphi},~ \widetilde\varphi,~ \varphi^L,~ \varphi^{\#}$}

\end{picture}
\end{math} \\ \\

and

\begin{math}
\setlength{\unitlength}{0.1cm}
\thicklines
\begin{picture}(50,31)

\put(0,12){$(4.2)$}
\put(23,25){$X \ni x$}
\put(62,29){$\varphi$}
\put(38,26){\vector(1,0){52}}
\put(38,25){\line(0,1){2}}
\put(95,25){$\varphi ( x ) \in Y$}
\put(33,22){\vector(0,-1){16}}
\put(32,22){\line(1,0){2}}
\put(17,0){$X^{\#} \ni \,\, < x \,]$}
\put(40,1){\vector(1,0){49}}
\put(40,0){\line(0,1){2.5}}
\put(93,0){$< \varphi ( x ) \,] \in Y^{\#}$}
\put(100,22){\vector(0,-1){16}}
\put(99,22){\line(1,0){2}}
\put(52,-4){$\overline{\varphi},~ \widetilde\varphi,~ \varphi^L,~ \varphi^{\#}$}

\end{picture}
\end{math} \\ \\

for every {\it cofinal} mapping $L$ in (2.4). \\

{\bf Proof} \\

It follows easily from the results in section 3. \\

{\bf Remark 4.1.} \\

The extension in (4.1) does in fact {\it not} need the partial order on $X$, and it comes down
to the extension in (A.28). \\

The extension in (4.2) comes down to the extension (A.29). \\

It follows that the extensions $\overline{\varphi},~ \widetilde\varphi$ and $\varphi^L$,
although not necessarily identical in general, do nevertheless reduce to $\varphi^{\#}$, in
the case of the diagrams (4.1) and (4.2). \\ \\

{\bf Appendix} \\

We shortly present several notions and results used above. A related full presentation can be
found in [3, Appendix, pp. 391-420]. \\

Let $(X,\leq)$ be a nonvoid poset without minimum or maximum. For $a \in X$ we denote \\

(A.1) $~~~ < a ] = \{ x \in X ~|~ x \leq a \},~~~ [ a > = \{ x \in X ~|~ x \geq a \} $ \\

We define the mappings \\

(A.2) $~~~ X ~\supseteq~ A \longmapsto A^u = \bigcap_{a \in A}~ [ a > ~\subseteq~ X $ \\

(A.3) $~~~ X ~\supseteq~ A \longmapsto A^l = \bigcap_{a \in A} < a ] ~\subseteq~ X $ \\

then for $A \subseteq X$ we have \\

(A.4) $~~~ A^u = X \Longleftrightarrow A^l = X \Longleftrightarrow A = \phi $ \\

(A.5) $~~~ A^u = \phi \Longleftrightarrow A ~\mbox{unbounded from above} $ \\

(A.6) $~~~ A^l = \phi \Longleftrightarrow A ~\mbox{unbounded from below} $ \\

{\bf Definition A.1.} \\

We call $A \subseteq X$ a {\it cut}, if and only if \\

(A.7) $~~~A^{ul} = A $ \\

and denote \\

(A.8) $~~~ X^{\#} = \{ A \subseteq X ~|~ A ~\mbox{is a cut} \} \subseteq {\cal P} ( X) $ \\

\hfill $\Box$ \\

Clearly, (A.4) - (A.6) imply \\

(A.9) $~~~ \phi,~ X \in X^{\#} $ \\

therefore \\

(A.10) $~~~ X^{\#} \neq \phi $ \\

Given $A, B \subseteq X$, we have \\

(A.11) $~~~ A \subseteq B \Longrightarrow A^u \supseteq B^u,~ A^l \supseteq B^l $ \\

(A.12) $~~~ A \subseteq A^{ul},~~~ A \subseteq A^{lu} $ \\

(A.13) $~~~ A^{ulu} = A^u,~~~ A^{lul} = A^l $ \\

Consequently \\

(A.14) $ \begin{array}{l}
            \forall~~ A \subseteq X ~: \\ \\
            ~~~~*)~~ A^{ul} \in X^{\#} \\ \\
            \begin{array}{l}
                    ~**)~~   \forall~~ B \in X^{\#} ~: \\ \\
                       ~~~~~~~~~~~~ A \subseteq B \Longrightarrow A^{ul} \subseteq B \\ \\
                       ~~~~~~~~~~~~ B \subseteq A \Longrightarrow B \subseteq A^{ul}
                     \end{array}
          \end{array} $ \\

therefore \\

(A.15) $~~~ X^{\#} = \{ A^{ul} ~|~ A \subseteq X \} $ \\

Given $x \in X$, we have \\

(A.16) $~~~ \{ x \}^u = [ x >,~~~ \{ x \}^l = < x ],~~~ [ x >^l = < x ],~~~
                                                < x ]^u = [ x > $ \\

(A.17) $~~~ \{ x \}^{ul} = < x],~~~ \{ x \}^{lu} = [ x > $ \\

We denote for short \\

$ \{ x \}^u = x^u,~~ \{ x \}^l = x^l,~~ \{ x \}^{ul} = x^{ul},~~
                                         \{ x \}^{lu} = x^{lu},~.~.~.~ $ \\

Given $A \in X^{\#}$, we have \\

(A.18) $~~~ \phi \neq A \neq X \Longleftrightarrow
                   \left (~ \begin{array}{l}
                               \exists~~ a, b \in X ~: \\ \\
                               ~~~ < a ] ~\subseteq~ A ~\subseteq~ < b ]
                            \end{array} ~~\right ) $ \\

We shall use the {\it embedding} \\

(A.19) $~~~ X \ni x ~\stackrel{\varphi}\longmapsto~ x^{ul} = x^l = < x ] \in X^{\#} $ \\

We define on $X^{\#}$ the partial order \\

(A.20) $~~~ A \leq B \Longleftrightarrow A \subseteq B $ \\ \\

{\bf Definition 2.1.} \\

Given two posets $( X, \leq ),~ ( Y, \leq )$ and a mapping $\varphi : X \longrightarrow Y$. We
call $\varphi$ an {\it order isomorphic embedding}, or in short, OIE, if and only if it is
injective, and furthermore, for $a, b \in X$ we have \\

$~~~ a ~\leq~ b ~~~\Longleftrightarrow~~~ \varphi ( a ) ~\leq~ \varphi ( b ) $ \\

An OIE $\varphi$ is an {\it order isomorphism}, or in short, OI, if and only if it is
bijective, which in this case is equivalent with being surjective.

\hfill $\Box$ \\

The main result concerning order completion is given in, [2] : \\ \\

{\bf Theorem ( H M MacNeille, 1937 )} \\

1)~ The poset $( X^{\#}, \leq )$ is order complete. \\

2)~ The embedding $X \stackrel{\varphi}\longrightarrow X^{\#}$ in (A.19) preserves infima and
suprema, and it is an order isomorphic embedding, or OIE. \\

3)~ For $A \in X^{\#}$, we have the order density property of $X$ in $X^{\#}$, namely \\

(A.21) $~~~ \begin{array}{l}
                A ~=~ \sup_{X^{\#}}~ \{ x^l ~|~ x \in X,~~ x^l \subseteq A \} ~=~ \\ \\
                  ~~~=~ \inf_{X^{\#}}~ \{ x^l ~|~ x \in X,~~ A \subseteq x^l \}
             \end{array} $ \\

\hfill $\Box$ \\

For $A \subseteq X$, we have \\

(A.22) $~~~ A^{ul} = \sup_{X^{\#}}~ \{ x^l ~|~ x \in A \} $ \\

Given $A_i \in X^{\#}$, with $i \in I$, we have with the partial order in $X^{\#}$ the
relations \\

(A.23) $~~~ \sup_{i \in I}~ A_i ~=~ \inf~ \{ A \in X^{\#} ~|~ \bigcup_{i \in I} A_i
                                     \subseteq A \} ~=~ ( \bigcup_{i \in I} A_i )^{ul} $ \\

(A.24) $~~~ \begin{array}{l}
                \inf_{i \in I}~ A_i ~=~ \sup~ \{ A \in X^{\#} ~|~ A \subseteq
                  \bigcap_{i \in I} A_i \} ~=~ ( \bigcap_{i \in I} A_i )^{ul} ~=~  \\ \\
               ~~~~~~~~~~~~~~=~ \bigcap_{i \in I} A_i
             \end{array} $ \\ \\

{\bf Extending mappings to order completions} \\

Let $(X,\leq),~(Y,\leq)$ be two posets without minimum or maximum, and let \\

(A.25) $~~~ \varphi : X \longrightarrow Y $ \\

be any mapping. Our interest is to obtain an extension \\

$~~~~~~ \varphi^{\#} : X^{\#} \longrightarrow Y^{\#} $ \\

For that, we first extend $\varphi$ to a {\it larger} domain, as follows \\

(A.26) $~~~ \varphi^{\#} : {\cal P} ( X ) \longrightarrow Y^{\#} $ \\

where for $A \subseteq X$ we define \\

(A.27) $~~~ \varphi^{\#} ( A ) ~=~ ( \varphi ( A ) )^{ul} ~=~
                \sup_{Y^{\#}}~ \{ < \varphi ( x )\, ] ~|~ x \in A \} $ \\

and for any mapping in (A.25), we obtain the commutative diagram \\

\begin{math}
\setlength{\unitlength}{0.1cm}
\thicklines
\begin{picture}(50,31)

\put(0,12){$(A.28)$}
\put(15,25){$X \ni x$}
\put(50,28){$\varphi$}
\put(28,26){\vector(1,0){48}}
\put(28,25){\line(0,1){2}}
\put(80,25){$\varphi ( x ) \in Y$}
\put(24.5,22){\vector(0,-1){17}}
\put(23.5,22){\line(1,0){2}}
\put(6,0){${\cal P} ( X ) \ni \{ x \}$}
\put(29,1){\vector(1,0){47}}
\put(29,0){\line(0,1){2}}
\put(80,0){$\varphi^{\#}( x ) \,=~ < \varphi ( x )\, ] \in Y^{\#}$}
\put(83,22){\vector(0,-1){17}}
\put(82,22){\line(1,0){2}}
\put(50,-3.5){$\varphi^{\#}$}

\end{picture}
\end{math} \\ \\

{\bf Proposition A.1.} \\

1)~ The mapping $\varphi^{\#} : {\cal P} ( X ) \longrightarrow Y^{\#}$ in (A.36) is
increasing, if on ${\cal P} ( X )$ we take the partial order defined by the usual inclusion
"$\subseteq$". \\

2)~ If the mapping $\varphi : X \longrightarrow Y$ in (A.35) is increasing, then the mapping
$\varphi^{\#} : {\cal P} ( X ) \longrightarrow Y^{\#}$ in (A.36) is an extension of it to
$X^{\#}$, namely, we have the commutative diagram \\

\begin{math}
\setlength{\unitlength}{0.1cm}
\thicklines
\begin{picture}(50,31)

\put(0,12){$(A.29)$}
\put(15,25){$X \ni x$}
\put(50,28){$\varphi$}
\put(28,26){\vector(1,0){48}}
\put(28,25){\line(0,1){2}}
\put(80,25){$\varphi ( x ) \in Y$}
\put(24.5,22){\vector(0,-1){17}}
\put(23.7,22){\line(1,0){2}}
\put(8,0){$X^{\#} \ni\, < x ]$}
\put(29,1){\vector(1,0){44}}
\put(29,0){\line(0,1){2}}
\put(77,0){$\varphi^{\#}( < x \,]\, ) \,=~ < \varphi ( x )\, ] \in Y^{\#}$}
\put(83,22){\vector(0,-1){17}}
\put(82,22){\line(1,0){2}}
\put(50,-4){$\varphi^{\#}$}

\end{picture}
\end{math} \\ \\

3)~ If the mapping $\varphi : X \longrightarrow Y$ in (A.25) is an OIE, then the mapping
$\varphi^{\#} : {\cal P} ( X ) \longrightarrow Y^{\#}$ in (A.26) when restricted to
$X^{\#}$, that is \\

(A.30) $~~~ \varphi^{\#} : X^{\#} \longrightarrow Y^{\#} $ \\

as in (A.29), is also an OIE. \\ \\

{\bf Lemma A.1.} \\

Let in general $\mu : M \longrightarrow N$ be an increasing mapping between two order
complete posets, then for nonvoid $E \subseteq M$ we have \\

(A.31) $~~~\mu (\, \inf_M\, E \,) ~\leq~ \inf_N\, \mu ( E )~\leq~
\sup_N\, \mu ( E )
                                                     ~\leq~ \mu (\, \sup_M\, E \,) $ \\

{\bf Proof} \\

Indeed, let $a = \inf_M\, E \in M$. Then $a \leq b$, with $b \in E$. Hence $\mu ( a ) \leq
\mu ( b )$, with $b \in E$. Thus $\mu ( a ) \leq \inf_N\, \mu ( E )$, and the first inequality
is proved. \\
The last inequality is obtained in a similar manner, while the middle inequality is trivial.

\hfill $\Box$ \\

\end{document}